# COMPUTING THE FROBENIUS-SCHUR INDICATOR FOR ABELIAN EXTENSIONS OF HOPF ALGEBRAS

Y. KASHINA, G. MASON, AND S. MONTGOMERY

## 1. INTRODUCTION

Let $H$ be a finite-dimensional semisimple Hopf algebra. Recently it was shown in [LM] that a version of the Frobenius-Schur theorem holds for Hopf algebras, and thus that the Schur indicator $\nu(\chi)$ of the character $\chi$ of a simple $H$-module is well-defined; this fact for the special case of Kac algebras was shown in [FGSV]. In this paper we show that for an important class of non-trivial Hopf algebras, $\nu(\chi)$ is a computable invariant. The Hopf algebras we consider are all abelian extensions; as a special case, they include the Drinfel'd double of a group algebra.

In addition to finding a general formula for the indicator, we also study when it is always positive. In particular we prove that the indicator is always positive for the Drinfeld double of the symmetric group, generalizing the classical result for the symmetric group itself. As a first step in proving this, we show that the indicator can be computed by means of a "local indicator".

Finally we show that work of the first author on the classification of Hopf algebras of dimension 16 can be somewhat shortened using indicators rather than $K_0$.

It is likely that the indicator will be useful in other problems on the classification of semisimple Hopf algebras. Moreover, Schur indicators play a role in various aspects of conformal field theory; see work of Bantay [B1] [B2].

We first introduce some notation. Throughout, $H$ will be a finite-dimensional Hopf algebra over an algebraically closed field $k$ of characteristic not 2, with comultiplication $\Delta : H \longrightarrow H \otimes H$, via $h \mapsto \sum h_1 \otimes h_2$, counit $\varepsilon : H \longrightarrow k$, and antipode $S$. We let $\mathcal{G}(H)$ denote the group of group-like elements of $H$. For a general reference on Hopf algebras, see [Mo].

We will assume throughout this paper that $H$ is semisimple, and if char $k$ is $p \neq 0$, that $H^*$ is also semisimple (in characteristic 0 this fact follows automatically by [LR]). In particular, by Maschke's theorem, $H$ has a (unique) integral $\Lambda$ such that $\varepsilon(\Lambda) = 1$.

The result on indicators we shall use is the following:

**Theorem 1.1.** [LM] *Let $H$ be a semisimple Hopf algebra over an algebraically closed field $k$. If $k$ has characteristic $p \neq 0$, assume in addition that $p \neq 2$ and that $H^*$ is semisimple. Let $\Lambda$ be an integral of $H$ with $\varepsilon(\Lambda) = 1$, and set $\Lambda^{[2]} := \sum_{(\Lambda)} \Lambda_1 \Lambda_2$. For a simple $H$-module $V$ with character $\chi$, define $\nu(\chi) := \chi(\Lambda^{[2]})$. Then the following properties hold:*

The first and third authors thank the Mathematical Sciences Research Institute, where some of this work was done, for support; also the third author was supported by NSF grant DMS 0100461.

(1) $\nu(\chi) = 0$, 1 or -1, for all such $\chi$.

(2) $\nu(\chi) \neq 0$ if and only if $V \cong V^*$. Moreover $\nu(\chi) = 1$ (respectively $-1$) if and only if $V$ admits a symmetric (resp. skew-symmetric) non-degenerate bilinear $H$-invariant form.

(3) Considering $S$ as an element of $End(H)$, $Trace(S) = \sum_\chi \nu(\chi)\chi(1_H)$, where the sum is over all simple characters $\chi$.

As for groups, we will call $\nu(\chi)$ the *Frobenius-Schur indicator* of $\chi$, or simply the *Schur indicator*. We frequently write $\nu(V)$ instead of $\nu(\chi)$. The theorem clearly specializes to the classical result for groups, noting that when $H = kG$ for a finite group $G$, $Sg = g^{-1}$, and thus the trace of $S$ is the number of involutions of $G$.

We note that the indicator may be viewed as a homomorphism of additive abelian groups

$$\nu : K_0(H) \longrightarrow \mathbb{Z}.$$

with $\nu$ taking the values 0, 1, or -1 when applied to a simple $H$-module $V$.

## 2. Abelian extensions

In order to describe the Hopf extensions in which we are interested, we first consider arbitrary extensions of finite-dimensional Hopf algebras. Thus

$$K \xhookrightarrow{i} H \xtwoheadrightarrow{\pi} F$$

where $K$, $H$, and $F$ are finite-dimensional Hopf algebras with $K$ normal in $H$ and $F = H/HK^+$. Since $H$ is finite-dimensional, it is known that $H = K\#_\sigma^\tau F$, with an action $\rightharpoonup : F \otimes K \to K$, a coaction $\rho : F \to F \otimes K$, a Hopf cocycle $\sigma : F \otimes F \to K$, and a dual cocycle $\tau : F \to K \otimes K$. See [A, 3.1.12] for details. As an algebra, $H = K\#_\sigma F$ is a crossed product of $F$ over $K$; thus writing the basis elements of $H$ as $w\#f = w\overline{f}$, where $w \in K$ and $f \in F$, the multiplication in $H$ is given by

$$(w\overline{f})(l\overline{g}) = \sum_{f,l} w(f_1 \rightharpoonup l)\sigma(f_2, g_1)\overline{f_3 g_2}.$$

The comultiplication and antipode in $H$ will be discussed below.

Note that the dual Hopf algebra $H^*$ is of the form $H^* = F^*\#_{\tau^*}^{\sigma^*}K^*$, where now $F^*$ is normal in $H^*$, $\tau^* : K^* \otimes K^* \to F^*$ is a Hopf cocycle on $K^*$, and $\sigma^*$ is a dual cocycle.

Since homomorphic images and Hopf subalgebras of semisimple Hopf algebras are semisimple (see [Mo]), it follows that $K, F, K^*$, and $F^*$ are also semisimple, by our assumption on $H$ and $H^*$.

The extension is called *abelian* if $K$ is commutative and $F$ is cocommutative; since $k$ is algebraically closed, it follows by [Mo, 2.3.1] that $K \cong (kG)^*$ and $F = H/HK^+ \cong kL$, for two groups $G$ and $L$. Thus we may assume that our extension is of the form

$$(2.1) \qquad (kG)^* \xhookrightarrow{i} H \xtwoheadrightarrow{\pi} kL,$$

where $H = K\#_\sigma^\tau kL$ as above. Moreover $\sigma$ and $\tau^*$ are simply group 2-cocycles twisted by the action; that is

(2.2) $$(z \rightharpoonup \sigma(x,y))\sigma(z, xy) = \sigma(z, x)\sigma(zx, y),$$

for $x, y, z \in L$, and similarly for $\tau^*$. We also assume that our cocycles are normalized; that is, $\sigma(1, x) = \sigma(x, 1) = 1$ for all $x \in G$. Also, since $H$ is a Hopf algebra, $\varepsilon(\sigma(x, y)) = \sigma(\varepsilon(x), \varepsilon(y)) = \sigma(1, 1) = 1$ for all $x, y \in G$ (again, see [A]).

We remark that by Maschke's Theorem, an abelian extension $H$ as in (2.1) will always be both semisimple and cosemisimple in characteristic 0, and will be semisimple and cosemisimple in characteristic $p > 0 \iff p$ does not divide $|G||L| = \dim H$.

Let $\{p_g | g \in G\}$ be the dual basis for $(kG)^*$. The action $\rightharpoonup$ of $L$ on $K$ induces an action of $L$ on $K^* = kG$ via $(l \rightharpoonup f)(k) := f(Sl \rightharpoonup k)$. Since $K$ is commutative and $kL$ is cocommutative, $K$ is a $kL$-module algebra, and thus $L$ acts as automorphisms of $K$. Thus $L$ permutes the orthogonal idempotents $\{p_g\}$; it follows that the action of $L$ on $kG$ is in fact an action of $L$ on $G$ itself, which we also denote by $\rightharpoonup$. Then the action of $L$ on $(kG)^*$ is given by

(2.3) $$x \rightharpoonup p_g = p_{x \rightharpoonup g}.$$

In order to compute with the cocycle $\sigma$, we write it in terms of the basis in $(kG)^*$. That is,

(2.4) $$\sigma(x, y) = \sum_{g \in G} \sigma_g(x, y) p_g$$

where $\sigma_g(x, y) \in k$. It is easy to see that $\sigma_1(x, y) = 1$, and that if $\sigma$ is trivial, then all $\sigma_g(x, y) = 1$. As a consequence of (2.2), we have

(2.5) $$\sigma_{z^{-1} \rightharpoonup g}(x, y)\sigma_g(z, xy) = \sigma_g(z, x)\sigma_g(zx, y).$$

Multiplication in $H$ can now be written as

(2.6) $$p_k \overline{z} p_h \overline{y} = p_k p_{z \rightharpoonup h} \sigma(z, y)\overline{zy} = \delta_{k, z \rightharpoonup h} \sigma_k(z, y) p_k \overline{zy}$$

where $h, k \in G$, $y, z \in L$. In particular, $\overline{z} p_h = p_{z \rightharpoonup h} \overline{z}$.

The comultiplication in an abelian extension is rather complicated. Thus most of our results are for cocentral abelian extensions. An extension is called *cocentral abelian* if it is abelian and in addition $F^* \subset Z(H^*)$. It follows that in $H^*$, the action of $K^*$ on $F^*$ is trivial, and thus in $H = K\#_\sigma^\tau kL$, the coaction $kL \to kL \otimes K$ is trivial. The comultiplication in such an $H$ is given by:

$$\Delta(p_g \overline{x}) = \Delta(p_g)\Delta(\overline{x}) = (\sum_{h \in G} p_h \otimes p_{h^{-1}g})\tau(x)(\overline{x} \otimes \overline{x}).$$

The counit of $H$ is given simply by $\varepsilon(p_g \overline{x}) = \delta_{g,1}$

As we did for $\sigma$ in 2.4, we may write $\tau$ in terms of the basis elements of $(kG)^*$. That is, $\tau(x) = \sum_{g,h \in G} \tau_{g,h}(x) p_g \otimes p_h$ for $\tau_{g,h}(x) \in k$. We may then write the comultiplication as

(2.7) $$\Delta(p_g \overline{x}) = \sum_{h \in G} \tau_{h, h^{-1}g}(x) p_h \overline{x} \otimes p_{h^{-1}g} \overline{x}.$$

Some further properties of $\tau$ are given in Lemma 4.5.

A general formula for the antipode is also given in [A]; it is rather complicated in general. In the case $H$ is cocentral, this specializes to

$$(2.8) \qquad S(p_g \overline{x}) = \sigma^{-1}_{x^{-1} \rightharpoonup g^{-1}}(x^{-1}, x) \tau^{-1}_{g^{-1}, g}(x) p_{x^{-1} \rightharpoonup g^{-1}} \overline{x^{-1}}$$

where we have written $\sigma^{-1}$ and $\tau^{-1}$ in terms of the basis for $(kG)^*$ as we did for $\sigma$ and $\tau$. As a special case, we note that when $\sigma$ and $\tau$ are trivial, we simply have $S(p_g x) = p_{x^{-1} \rightharpoonup g^{-1}} x^{-1}$.

As remarked earlier, the Drinfeld double is an example of the extensions we study.

**Example 2.9.** *The Drinfeld double of a group algebra*

The Drinfeld double $H = D(G)$ of a group $G$ is just a cocentral abelian extension as above with $L = G$ such that $G$ acts on itself by conjugation and with trivial cocycle and dual cocycle. Thus $x \rightharpoonup g = xgx^{-1}$ and $x \rightharpoonup p_g = p_{xgx^{-1}}$, for $x, g \in G$.

Writing the basis elements of $D(G)$ as $p_g \bowtie x = p_g \# x$, for $g, x \in G$, multiplication is given by $(p_k \bowtie z)(p_h \bowtie y) = \delta_{k, zhz^{-1}} p_k \bowtie zy$, and comultiplication is given by $\Delta(p_g \bowtie x) = \sum_{h \in G} (p_h \bowtie x) \otimes (p_{h^{-1}g} \bowtie x)$, as in (2.7).

The antipode is given by $S(p_g \bowtie x) = p_{x^{-1}g^{-1}x} \bowtie x^{-1}$.

We require some additional definitions concerning the action of $L$ on $G$. For $g, h \in G$, let $L_{g,h} = \{y \in L | y \rightharpoonup g = h\}$; then $L_{g,g} = L_g$, the *stabilizer* of $g$ in $L$. Later on we will need an extension of $L_g$; that is, let $\tilde{L}_g := L_g \cup L_{g,g^{-1}}$, the *extended stabilizer* of $g$ in $L$. Note that $L_{g,g^{-1}} = L_{g^{-1},g}$ and that $(L_{g,g^{-1}})^2 \subset L_g$. Thus $\tilde{L}_g$ is a subgroup of $L$ with $[\tilde{L}_g : L_g] \leq 2$.

We also say that $g \in G$ is *L-real* if $L_{g,g^{-1}} \neq \emptyset$, and *L-non-real* if $L_{g,g^{-1}} = \emptyset$. Finally we let $\mathcal{O}(g) = \{y \rightharpoonup g | y \in L\}$ denote the *orbit* of $g$ under the action of $L$ and let $T_g$ denote a complete set of left coset representatives of $L_g$ in $L$.

When $L = G$ acts on itself by conjugation, as in the case of $H = D(G)$, then $L_g = C_g$, the centralizer of $g$ in $G$, and $\mathcal{O}(g)$ is simply the conjugacy class of $g$. We also write $C_{g,g^{-1}}$ for $L_{g,g^{-1}}$. Then $x$ is *real* if it is inverted under conjugation by some element of $G$, and *non-real* otherwise.

## 3. Modules over crossed products

In this section, we do not need the fact that $H$ is a Hopf algebra. Thus we only assume that $H = (kG)^* \#_\sigma kL$, a crossed product of the group $L$ over the dual group algebra $K = (kG)^*$, where $\sigma$ is an ordinary 2-cocycle from $L$ to the commutative ring $K$, with an $L$-action, and $L$ acts on $K$ via a given action on $G$. That is, letting $p_g$, for $g \in G$, be the usual basis of orthogonal idempotents for $K$ as in the previous section, we assume that $y \rightharpoonup p_g = p_{y \rightharpoonup g}$, for any $y \in L$, as in (2.3), that $\sigma$ satisfies (2.2), and that the multiplication in $H$ is given by (2.6).

We will show that all simple left $H$-modules may be described as induced modules from certain twisted group algebras of the stabilizers $L_g$. Our approach here follows the general outline of the work on more general crossed products $B = A \#_\sigma kL$, for any semisimple algebra $A$, in [MoW] and [W]. However here we get more explicit information more easily, because of the special nature of our algebra $A = (kG)^*$ and of the $L$-action on $A$.

First note that, writing the cocycle $\sigma$ using the basis $\{p_g\}$ as in (2.4), it follows from (2.5) that $\sigma_g$ is an ordinary 2-cocycle when restricted to the stabilizer $L_g$. Thus we may form the twisted group algebra $k_{\sigma_g} L_g$. A basis of $k_{\sigma_g} L_g$ is given by $\{\tilde{x} \mid x \in L_g\}$, and multiplication is given by

$$\tilde{x}\tilde{y} = \sigma_g(x,y)\widetilde{xy}. \tag{3.1}$$

**Lemma 3.2.** *Let $H = (kG)^* \#_\sigma kL$ be a crossed product as above, and fix an element $g \in G$. Let $k_{\sigma_g} L_g$ be the twisted group algebra defined above. Then*

*(1) $k_{\sigma_g} L_g$ is isomorphic as a $k$-algebra to the subalgebra $p_g((kG)^* \#_\sigma kL_g)p_g = p_g \#_\sigma kL_g$ of $H$.*

*(2) Let $V$ be a left $k_{\sigma_g} L_g$-module, and let $V' = p_g \otimes V$. Then $V'$ becomes a left $(kG)^* \#_\sigma kL_g$-module by defining, for each $h \in G$, $x \in L_g$, $v \in V$,*

$$(p_h \overline{x}) \cdot (p_g \otimes v) := \delta_{h,g}(p_g \otimes \tilde{x} \cdot v).$$

*(3) Let $\hat{V} := H \otimes_{(kG)^* \#_{\sigma_g} kL_g} V' = H \otimes_{(kG)^* \#_{\sigma_g} L_g} (p_g \otimes V)$. Then $\hat{V}$ is a left $H$-module via left multiplication, the usual induced $H$-module structure. Moreover we may write $\hat{V} = \sum_{z \in T_g} \overline{z} \otimes (p_g \otimes V)$.*

*Combining these actions, the $H$-action on $\hat{V}$ is given by the following formula, where $h \in G$, $x, y \in L$, and $v \in V$:*

$$(p_h \overline{x}) \cdot [\overline{y} \otimes (p_g \otimes v)] = \delta_{h,xy \to g} \sigma_{xy \to g}(x,y) \overline{xy} \otimes (p_g \otimes v).$$

*Proof.* (1) Since for $x \in L_g$, $p_g \overline{x} = \overline{x} p_g$, it is clear that $p_g((kG)^* \#_\sigma kL_g)p_g = p_g \#_\sigma kL_g$. The isomorphism $\Phi: k_{\sigma_g} L_g \to p_g \#_\sigma kL_g$ is given by $\Phi(\tilde{x}) = p_g \overline{x}$, since for $x, y \in L_g$,

$$\Phi(\tilde{x})\Phi(\tilde{y}) = (p_g \overline{x})(p_g \overline{y}) = p_g \sigma_g(x,y)\overline{xy} = \Phi(\tilde{x}\tilde{y})$$

(2) The fact that $V'$ is a $(kG)^* \#_\sigma kL_g$-module follows from the algebra isomorphism in (1).

(3) Clearly $\hat{V}$ is a left $H$-module. Since $L = \cup_{z \in T_g} L_g$, $H$ is a free right $(kG)^* \#_\sigma L_g$-module with basis $T_g$, and so $\hat{V} = \sum_{z \in T_g} \overline{z} \otimes (p_g \otimes V)$. □

For the rest of the paper, we will fix our notation as in Lemma 3.2; that is, for a given $g \in G$, $V = V_g$ will denote a $k_{\sigma_g} L_g$-module, $V' = p_g \otimes V$ will denote the corresponding $(kG)^* \#_\sigma kL_g$-module, and $\hat{V} = H \otimes_{(kG)^* \#_\sigma kL_g} V'$ the induced $H$-module.

**Theorem 3.3.** *Let $H = (kG)^* \#_\sigma kL$ be a crossed product as above. For each $L$-orbit $\mathcal{O}$ of $G$, fix an element $g \in \mathcal{O}$ and let $V$ be a left $k_{\sigma_g} L_g$-module. Let $\hat{V} = H \otimes_{(kG)^* \#_\sigma kL_g} (p_g \otimes V)$, the induced $H$-module as above. If $V$ is a simple $k_{\sigma_g} L_g$-module, then $\hat{V}$ is a simple $H$-module.*

*Conversely every simple left $H$-module is isomorphic to $\hat{V}$ for some simple module $V$ of $k_{\sigma_g} L_g$, where $g$ ranges over a choice of one element in each $L$-orbit $\mathcal{O}$ of $G$.*

*Proof.* From Lemma 3.2 we know that $\hat{V}$ is in fact an $H$-module. Now assume that $V$ is simple; we claim $\hat{V}$ is simple. Let $U$ be an $H$-submodule of $\hat{V}$ and let $0 \neq w \in U$. Then we may write $w = \sum_{x \in T_g} \overline{x} \otimes (p_g \otimes w_x) \in U$, where $w_z \neq 0$ for some $z \in T_g$. Consider $(p_g \overline{z^{-1}}) \cdot w \in U$:

$$\begin{aligned}
(p_g \overline{z^{-1}}) \cdot w &= (p_g \# z^{-1}) \cdot (\sum_{x \in T_g} \overline{x} \otimes (p_g \otimes w_x)) \\
&= \sum_{x \in T_g} \delta_{g, z^{-1}x \rightharpoonup g} \sigma_{z^{-1}x \rightharpoonup g}(z^{-1}, x)(\overline{z^{-1}x} \otimes (p_g \otimes w_x)) \\
&= \sigma_g(z^{-1}, z)(\overline{1} \otimes (p_g \otimes w_z)),
\end{aligned}$$

since $z^{-1}x \in L_g$ only if $z = x$.

Since $\sigma$ is convolution invertible, $\sigma_g(z^{-1}, z) \neq 0$ and therefore $\overline{1} \otimes p_g \otimes w_z \in U$. Since $V$ is a simple $k_{\sigma_g}L_g$-module, $k_{\sigma_g}L_g \cdot w_z = V$. Thus $\overline{1} \otimes (p_g \otimes V) = (1 \# kL_g) \cdot (\overline{1} \otimes (p_g \otimes w_z)) \subseteq U$. Therefore for any $x \in T_g$ and $v \in V$, $\overline{x} \otimes (p_g \otimes v) = \overline{x} \cdot (\overline{1} \otimes p_g \otimes v) \in U$ and thus $U = \hat{V}$.

If $W$ is a simple $k_{\sigma_h}L_h$-module, then it is easy to see that $\hat{V}_g$ and $\hat{W}_h$ are nonisomorphic unless $h \in \mathcal{O}(g)$. Moreover if $V_1$ and $V_2$ are non-isomorphic $k_{\sigma_g}L_g$-modules, then $\hat{V}_1$ and $\hat{V}_2$ are non-isomorphic as $H$-modules. We claim that all $\hat{V}_g$, where one element $g$ is chosen from each $L$-orbit of $G$, exhaust all possible simple $H$-modules.

To see this we use a counting argument. First fix $g \in G$. Then taking the sum over all nonisomorphic simple $k_{\sigma_g}L_g$-modules $V_g$, we obtain

$$\begin{aligned}
\sum_{V_g}(\dim \hat{V}_g)^2 &= \sum_{V_g}([L:L_g]\dim V_g)^2 = [L:L_g]^2 \sum_{V_g}(\dim V_g)^2 \\
&= [L:L_g]^2|L_g| = [L:L_g]|L| = |\mathcal{O}(g)||L|.
\end{aligned}$$

Now taking the sum over all distinct $L$-orbits $\mathcal{O}(g)$ we get

$$\sum_{\mathcal{O}(g)} \sum_{V_g}(\dim \hat{V}_g)^2 = |G||L| = \dim((kG)^* \#_\sigma^\tau kL).$$

Thus any simple $H$-module must be among the $\hat{V}_g$. □

It is clear from the theorem that the set of simple $H$-modules is a disjoint union over the distinct $L$-orbits of $G$ of those modules which are induced from $k_{\sigma_g}L_g$, with one $g$ chosen from each orbit. However a more precise statement can be made.

To see this, let $\mathcal{O}$ be an orbit of $L$ on $G$. Then we define

$$H(\mathcal{O}) := \sum_{g \in \mathcal{O}} p_g \# kL.$$

It follows from the multiplication in $H$ that the $H(\mathcal{O})$ are ideals of $H$ and that $H = \oplus_{\mathcal{O}} H(\mathcal{O})$. Note also that the underlying space of $\hat{V}$ is naturally an $H(\mathcal{O})$-module. We have:

**Corollary 3.4.** *Fix an element $g$ in the $L$-orbit $\mathcal{O}$ of $G$. The association $V \to \hat{V}$ in the theorem is the object map of a Morita equivalence of categories*

$$Mod_{k_{\sigma_g} L_g} \longrightarrow Mod_{H(\mathcal{O})}.$$

*This categorical equivalence induces an isomorphism $\phi$ between the corresponding Grothendieck groups $K_0(k_{\sigma_g} L_g)$ and $K_0(H(\mathcal{O}))$. Thus $K_0(H) = \oplus_\mathcal{O} K_0(H(\mathcal{O}))$.*

In the case when the cocycle $\sigma$ is trivial, Theorem 3.3 has a simpler form. In that case, the twisted group algebra is just $kL_g$ and we may identify $V'$ with $V$ and $\hat{V}$ with $Ind_{L_g}^L V = kL \otimes_{kL_g} V$. Then we have:

**Corollary 3.5.** *Let $H = (kG)^* \# kL$ be a smash product, where $L$ acts on $(kG)^*$ as in (2.3) and the multiplication is given by $xp_g = p_{x \to g} x$. Fix an element $g \in G$ and let $kL_g \cong 1 \# kL_g$ be the group algebra. Let $V$ be a left $kL_g$-module, and let $\hat{V} := Ind_{L_g}^L V = kL \otimes_{L_g} V$. Then $\hat{V}$ is a left $H$-module, via the action*

$$(p_h x) \cdot [y \otimes v] = \delta_{xy \to g, h} xy \otimes v.$$

*for $h \in G$, $x, y \in L$, and $v \in V$.*

*If $V$ is a simple left $kL_g$-module, then $\hat{V}$ is a simple left $H$-module. Conversely every simple left $H$-module is isomorphic to $\hat{V}$ for some simple module $V$ of $kL_g$, where $g$ ranges over a choice of one element in each $L$-orbit $\mathcal{O}$ of $G$.*

**Example 3.6.** Let $H = D(G)$ be the Drinfeld double of $G$, as in Example 2.9. As noted in the introduction, the orbits of $L$ acting on $G$ are the conjugacy classes of $G$, and the stabilizer $L_g = C_g$, the centralizer of $g$ in $G$. Thus as a special case of Theorem 3.3, or of Corollary 3.5, we see that the simple $D(G)$-modules arise by choosing one $g$ in each conjugacy class of $G$, and inducing the simple $C_g$-modules up to $D(G)$ as in 3.5. We thus recover the known facts about irreducible $D(G)$-modules from [Ma, Section 2], [DPR], and [AG, 3.1.1 and 3.1.2].

## 4. The Schur indicator for cocentral abelian extensions

From now on, we will assume that our Hopf algebra $H$ is a cocentral abelian extension. In this section we find a general formula for the Schur indicator for such extensions, although we get a more usable result when the cocycle $\sigma$ is trivial.

Before finding the Schur indicator, we must compute $\Lambda^{[2]} = \sum \Lambda_1 \Lambda_2$, where $\Lambda$ is an integral of $H$ with $\varepsilon(\Lambda) = 1$. For any (semisimple) extension, if $T$ and $t$ are integrals of $K$ and $F$ respectively, then it can be checked that $\Lambda = T \# t = T\bar{t}$ is an integral of $H$. In the abelian case, integrals for $K$ and $F$ are $T = p_1$ and $t = \frac{1}{|L|} \sum_{x \in L} x$ respectively, and so $H$ has integral $\Lambda = \frac{1}{|L|} \sum_{x \in L} p_1 \# x$ with $\varepsilon(\Lambda) = 1$.

**Lemma 4.1.** *Let $\Lambda = \frac{1}{|L|} \sum_{x \in L} p_1 \# x$. Then*

$$\Lambda^{[2]} = \frac{1}{|L|} \sum_{g \in G} \sum_{x \in L_{g,g^{-1}}} \tau_{g,g^{-1}}(x) \sigma_g(x,x) p_g \overline{x^2}.$$

*Proof.* We first compute $\Delta$ on $\Lambda$, using (2.7).

$$\Delta(\Lambda) = \frac{1}{|L|}\sum_{x\in L}\Delta(p_1\#x)$$

$$= \frac{1}{|L|}\sum_{g\in G}\sum_{x\in L}\tau_{g,g^{-1}}(x)p_g\overline{x}\otimes p_{g^{-1}}\overline{x}.$$

Thus

$$\Lambda^{[2]} = \frac{1}{|L|}\sum_{g\in G}\sum_{x\in L}\tau_{g,g^{-1}}(x)p_g\overline{x}p_{g^{-1}}\overline{x} = \frac{1}{|L|}\sum_{g\in G}\sum_{x\in L}\tau_{g,g^{-1}}(x)p_gp_{x\rightharpoonup g^{-1}}\overline{x}\,\overline{x}$$

$$= \frac{1}{|L|}\sum_{g\in G}\sum_{x\in L_{g^{-1},g}}\tau_{g,g^{-1}}(x)p_g\sigma(x,x)\overline{x^2}.$$

$\square$

We give an immediate application, which generalizes the classical fact in group theory that for a group of odd order, every non-trivial irreducible module has indicator 0.

**Corollary 4.2.** *Assume that $H$ is a cocentral abelian extension and that $\dim H$ is odd. Then for any simple character $\chi$ of $H$,*

$$\nu(\chi) = \chi(\Lambda^{[2]}) = \chi(\Lambda) = \begin{cases} 1 & \text{if } \chi \text{ is trivial} \\ 0 & \text{if } \chi \text{ is not trivial} \end{cases}$$

*In particular, a simple $H$-module is self-dual only if it is trivial.*

*Proof.* We apply Lemma 4.1. Since $\dim H$ is odd, $|L|$ and $|G|$ are odd. Thus $x \rightharpoonup g = g^{-1}$ implies $g = g^{-1}$ (otherwise $x$ has even order) and therefore $g = 1$ (otherwise $g$ has order 2). Thus

$$\Lambda^{[2]} = \frac{1}{|L|}\sum_{x\in L}\tau_{1,1}(x)\sigma_1(x,x)p_1\overline{x^2} = p_1\frac{1}{|L|}\sum_{x\in L}\overline{x^2} = T^{[2]}\#t^{[2]} = T\#t = \Lambda.$$

The result now follows since $\chi(\Lambda) = 0$ if $\chi$ is not trivial and $\chi(\Lambda) = 1$ if $\chi$ is trivial. $\square$

We note that in fact the conclusion of the corollary is true for any $H$ of odd dimension, even if $H$ is not an abelian extension [KSZ]. We give an alternate generalization of Corollary 4.2 in Theorem 4.13.

Note that the Corollary says that when $\dim H$ is odd, $\sum_\chi \nu(\chi)\chi(1_H) = 1$. The next proposition generalizes this fact.

**Proposition 4.3.** *Let $H$ be a cocentral abelian extension, as in (2.1). Then*

$$\sum_\chi \nu(\chi)\chi(1_H) = \sum_{\{(g,x)\in G\times L|x^2=1, x\rightharpoonup g=g^{-1}\}}\sigma_g^{-1}(x,x)\tau_{g,g^{-1}}^{-1}(x)$$

*where the sum is over all simple characters $\chi$ of $H$. In particular if both $\sigma$ and $\tau$ are trivial, then*

$$\sum_\chi \nu(\chi)\chi(1_H) = |\{(g,x)\in G\times L|x^2=1, x\rightharpoonup g=g^{-1}\}|.$$

*Proof.* This follows directly from the fact that $Tr(S) = \sum_\chi \nu(\chi)\chi(1_H)$ by Theorem 1.1(iii), the fact that the $p_g \overline{x}$, for all $g \in G, x \in L$, are a basis for $H$, and the formula for the antipode (2.8). $\square$

In the special case when $H = D(G)$, as in Example 2.9, the sum in Proposition 4.3 counts the number of maps of (a certain quotient of) the fundamental group of the Klein bottle into $G$; see [B2].

The next result is our most general, when both $\sigma$ and $\tau$ are non-trivial.

**Theorem 4.4.** *Let $H$ be cocentral abelian, as in (2.1). For $g \in G$, let $V$ be a simple $k_{\sigma_g} L_g$-module and let $\hat{V}$ be the corresponding simple $H$-module. Let $\chi$ denote the character of $V$ and let $\hat{\chi}$ be the character of $\hat{V}$. Then*

$$\nu(\hat{\chi}) = \frac{1}{|L|} \sum_{x \in T_g, z \in L_{g,g^{-1}}} \frac{\tau_{x \rightharpoonup g, x \rightharpoonup g^{-1}}(xzx^{-1})\sigma_{x \rightharpoonup g}(xzx^{-1}, xzx^{-1})\sigma_{x \rightharpoonup g}(xz^2x^{-1}, x)}{\sigma_{x \rightharpoonup g}(x, z^2)} \chi(z^2)$$

*Proof.* Let $h \in G$, and assume $y \in L_{h,h^{-1}} = L_{h^{-1},h}$. Then $y^2 \in L_{h,h^{-1}} L_{h^{-1},h} \subseteq L_h$ and so for any $x \in L$, $x^{-1}y^2 x \in L_g \iff x \in L_{g,h}$. Then for $v \in V$, $x \in T_g$, using Lemma 3.2,

$$\begin{aligned}
p_h \overline{y^2} \cdot (\overline{x} \otimes (p_g \otimes v)) &= \delta_{h, y^2 x \rightharpoonup g} \sigma_{y^2 x \rightharpoonup g}(y^2, x)(\overline{y^2 x} \otimes (p_g \otimes v)) \\
&= \delta_{h, x \rightharpoonup g} \sigma_{x \rightharpoonup g}(y^2, x)(\overline{y^2 x} \otimes (p_g \otimes v)) \\
&= \delta_{h, x \rightharpoonup g} \frac{\sigma_{x \rightharpoonup g}(y^2, x)}{\sigma_{x \rightharpoonup g}(x, x^{-1}y^2 x)} (\overline{x} \otimes (p_g \otimes \widetilde{x^{-1}y^2 x} \cdot v))
\end{aligned}$$

It follows that if $L_{g,h} = \emptyset$, then $\hat{\chi}(p_{x \rightharpoonup g} \overline{y^2}) = 0$. If $L_{g,h} \neq \emptyset$, then $L_{g,h} = xL_g$ for some $x \in L_{g,h}$, and thus

$$\hat{\chi}(p_h \overline{y^2}) = \hat{\chi}(p_{x \rightharpoonup g} \overline{y^2}) = \frac{\sigma_{x \rightharpoonup g}(y^2, x)}{\sigma_{x \rightharpoonup g}(x, x^{-1}y^2 x)} \chi(x^{-1}y^2 x).$$

It then follows from Lemma 4.1 that

$$\hat{\chi}(\Lambda^{[2]}) = \frac{1}{|L|} \sum_{h \in G} \sum_{y \in L_{h,h^{-1}}} \tau_{h,h^{-1}}(y)\sigma_h(y,y)\hat{\chi}(p_h \overline{y^2})$$

$$= \frac{1}{|L|} \sum_{x \in T_g, y \in xL_{g,g^{-1}}x^{-1}} \tau_{x \rightharpoonup g, x \rightharpoonup g^{-1}}(y)\sigma_{x \rightharpoonup g}(y,y)\hat{\chi}(p_{x \rightharpoonup g}\overline{y^2})$$

$$= \frac{1}{|L|} \sum_{x \in T_g, y \in xL_{g,g^{-1}}x^{-1}} \frac{\tau_{x \rightharpoonup g, x \rightharpoonup g^{-1}}(y)\sigma_{x \rightharpoonup g}(y,y)\sigma_{x \rightharpoonup g}(y^2, x)}{\sigma_{x \rightharpoonup g}(x, x^{-1}y^2 x)} \chi(x^{-1}y^2 x)$$

$$= \frac{1}{|L|} \sum_{x \in T_g, z \in L_{g,g^{-1}}} \frac{\tau_{x \rightharpoonup g, x \rightharpoonup g^{-1}}(xzx^{-1})\sigma_{x \rightharpoonup g}(xzx^{-1}, xzx^{-1})\sigma_{x \rightharpoonup g}(xz^2x^{-1}, x)}{\sigma_{x \rightharpoonup g}(x, z^2)} \chi(z^2)$$

$\square$

From now on we will usually assume that $\sigma$ is trivial, to simplify the situation. We require some further properties of $\tau$.

**Lemma 4.5.** *As above, let $\tau(x) = \sum_{g,h \in G} \tau_{g,h}(x) p_g \otimes p_h$, where $\tau_{g,h}(x) \in k$. Then the following properties are satisfied, for $x, y, z \in L$ and $g, h \in G$:*

$$\tau_{g,h}(x)\tau_{gh,k}(x) = \tau_{g,hk}(x)\tau_{h,k}(x) \tag{4.6}$$

$$\tau_{g,1}(x) = \tau_{1,g}(x) = 1 \tag{4.7}$$

$$\tau_{g,h}(x)\tau_{x^{-1} \rightharpoonup g, x^{-1} \rightharpoonup h}(y)\sigma_g(x,y)\sigma_h(x,y) = \tau_{g,h}(xy)\sigma_{gh}(x,y) \tag{4.8}$$

$$\tau_{g,g^{-1}}(x) = \tau_{g^{-1},g}(x) \tag{4.9}$$

*If also $\sigma$ is trivial and $z \in \tilde{L}_g$, then*

$$\tau_{x \rightharpoonup g, x \rightharpoonup g^{-1}}(xzx^{-1}) = \tau_{g,g^{-1}}(z) \tag{4.10}$$

$$\tau_{g,g^{-1}}(zx) = \tau_{g,g^{-1}}(z)\tau_{g,g^{-1}}(x) \tag{4.11}$$

*and thus $\tau_{g,g^{-1}} \in \mathcal{G}((k\tilde{L}_g)^*)$.*

*Proof.* The first formula is simply the dual cocycle condition, and the second is the fact that the dual cocycle is normalized. Property (4.8) follows from the fact that $\Delta$ is an algebra map, using (2.6) and (2.7). To see (4.9), we use (4.6) and (4.7):

$$\tau_{g,g^{-1}}(x) = \tau_{g,g^{-1}}(x)\tau_{gg^{-1},g}(x) = \tau_{g,g^{-1}g}(x)\tau_{g^{-1},g}(x) = \tau_{g^{-1},g}(x).$$

Now assume that $\sigma$ is trivial. By (4.8), it follows that

$$\tau_{g,g^{-1}}(z)\tau_{z^{-1} \rightharpoonup g, z^{-1} \rightharpoonup g^{-1}}(x) = \tau_{g,g^{-1}}(zx).$$

If $z \in L_g$ then $\tau_{z^{-1} \rightharpoonup g, z^{-1} \rightharpoonup g^{-1}}(x) = \tau_{g,g^{-1}}(x)$, and if $z \in \tilde{L}_g \setminus L_g$ then $\tau_{z^{-1} \rightharpoonup g, z^{-1} \rightharpoonup g^{-1}}(x) = \tau_{g^{-1},g}(x) = \tau_{g,g^{-1}}(x)$ by formula (4.9). Thus $\tau_{z^{-1} \rightharpoonup g, z^{-1} \rightharpoonup g^{-1}}(x) = \tau_{g,g^{-1}}(x)$ for any $z \in \tilde{L}_g$, proving (4.11).

Next, note that $\tau_{g,g^{-1}}(x^{-1})$ is a unit since by (4.8),

$$\tau_{g,g^{-1}}(x^{-1})\tau_{x \rightharpoonup g, x \rightharpoonup g^{-1}}(x) = \tau_{g,g^{-1}}(x^{-1}x) = 1$$

By formula (4.8) again,

$$\tau_{g,g^{-1}}(x^{-1})\tau_{x \rightharpoonup g, x \rightharpoonup g^{-1}}(xzx^{-1}) = \tau_{g,g^{-1}}(zx^{-1})$$

Therefore $\tau_{x \rightharpoonup g, x \rightharpoonup g^{-1}}(xzx^{-1}) = \tau_{g,g^{-1}}(z)$, using (4.11) and the fact that $\tau_{g,g^{-1}}(x^{-1})$ is a unit. This proves (4.10).

The map $\tau_{g,g^{-1}} : k\tilde{L}_g \to k$ is linear by the properties of the dual cocycle and so is in $(k\tilde{L}_g)^*$. Moreover, by (4.11), $\tau_{g,g^{-1}}$ is multiplicative and therefore $\tau_{g,g^{-1}} \in \mathcal{G}((k\tilde{L}_g)^*)$. □

**Corollary 4.12.** *Assume that $\sigma$ is trivial. Then*

$$\nu_H(\hat{\chi}) = \frac{1}{|L_g|} \sum_{z \in L_{g,g^{-1}}} \tau_{g,g^{-1}}(z)\chi(z^2)$$

*Proof.* Applying Theorem 4.4 and (4.10) we get:

$$\begin{aligned}
\nu_H(\hat\chi) &= \frac{1}{|L|} \sum_{x \in T_g, z \in L_{g,g^{-1}}} \tau_{x \to g, x \to g^{-1}}(xzx^{-1})\chi(z^2) \\
&= \frac{1}{|L|} \sum_{x \in T_g, z \in L_{g,g^{-1}}} \tau_{g,g^{-1}}(z)\chi(z^2) \\
&= \frac{1}{|L|} |T_g| \sum_{z \in L_{g,g^{-1}}} \tau_{g,g^{-1}}(z)\chi(z^2) \\
&= \frac{1}{|L_g|} \sum_{z \in L_{g,g^{-1}}} \tau_{g,g^{-1}}(z)\chi(z^2)
\end{aligned}$$

□

We can now give the generalization of Corollary 4.2 mentioned earlier. It demonstrates the importance of the involutions in $G$. Since we assume that $\sigma$ is trivial, we may use the simpler form of $\hat V$ in Corollary 3.5.

**Theorem 4.13.** *Let $H$ be cocentral abelian, and assume that $\sigma$ is trivial and $L$ is of odd order. Then for any $\hat\chi \in Irr(H)$, $\nu_H(\hat\chi)$ is nonnegative.*

*More precisely, let $\hat\chi$ be the character of $\hat V = Ind_{L_g}^L(V)$ for a simple $L_g$-module $V$, for $g \in G$, and let $\chi$ be the character of $V$. Then*

*(1) $\nu_H(\hat\chi) = 1 \iff g$ is an involution in $G$ and $\chi = \psi^{-1}$, where $\psi$ is the (unique) element in $\mathcal{G}((k\tilde L_g)^*)$ such that $\tau_{g,g^{-1}} = \psi^2$.*

*(2) If also $\tau$ is trivial, then $\nu_H(\hat\chi) = 1 \iff g$ is an involution in $G$ and $\chi$ is the trivial character for $L_g$.*

*Thus the number of simple characters $\hat\chi$ of $H$ which have $\nu_H(\hat\chi) = 1$ is exactly $n + 1$, where $n$ is the number of orbits of involutions in $G$.*

*Proof.* Let $g \in G$ and $\varphi = \tau_{g,g^{-1}}$. By Lemma 4.5, $\varphi \in \mathcal{G}((k\tilde L_g)^*)$, an abelian group of odd order. Therefore there exists a unique $\psi \in \mathcal{G}((k\tilde L_g)^*)$ such that $\tau_{g,g^{-1}} = \varphi = \psi^2$. For any simple character $\chi$ of $kL_g$, $\psi * \chi$ is another simple character of $kL_g$ of the same degree as $\chi$.

Let us now consider $\hat\chi$. Since $|L|$ is odd, we always have $L_g = \tilde L_g$ and there are only two possibilities:

(1) $o(g) \geq 3$ and $L_{g,g^{-1}} = \emptyset$. In this case $\nu_H(\hat\chi) = 0$ for any $\chi$.

(2) $g^2 = 1$ and $L_{g,g^{-1}} = L_g = \tilde L_g$.

Therefore by Corollary 4.12

$$\begin{aligned}
\nu_H(\hat{\chi}) &= \frac{1}{|L_g|} \sum_{z \in L_g} \tau_{g,g^{-1}}(z)\chi(z^2) = \frac{1}{|L_g|} \sum_{z \in L_g} \psi^2(z)\chi(z^2) \\
&= \frac{1}{|L_g|} \sum_{z \in L_g} \psi(z^2)\chi(z^2) = \frac{1}{|L_g|} \sum_{z \in L_g} (\psi * \chi)(z^2) \\
&= \nu(\psi * \chi) = \begin{cases} 1 & \text{if } \psi * \chi \text{ is trivial} \\ 0 & \text{otherwise} \end{cases}
\end{aligned}$$

Here we have used again the classical fact about indicators for characters of groups of odd order. Therefore if $g^2 = 1$ then $\nu_H(\hat{\chi}) = 0$ unless $\chi = \psi^{-1}$, in which case $\nu_H(\hat{\chi}) = 1$. □

We close this section with an example showing that a skew representation of $L_g$ can become real when induced up to $H$.

**Example 4.14.** Let $G = \langle a, b, c, d : a^4 = b^2 = c^4 = 1, d^2 = c^2, ba = a^{-1}b, ca = ac, da = ad, bc = c^{-1}b, dc = c^{-1}d, bd = d^{-1}b \rangle$ and let $H = D(G)$.

It is easy to check that $L_a = C_a = \langle a, c, d : a^4 = c^4 = 1, d^2 = c^2, dc = c^{-1}d \rangle$ and that $C_{a,a^{-1}} = bC_a$. Therefore for any irreducible $C_g$-module $V$ with character $\chi$,

$$\begin{aligned}
\nu_H(\hat{\chi}) &= \frac{1}{32} \sum_{i,j=0}^{3} \sum_{k=0}^{1} \chi((ba^i c^j d^k)^2) = \frac{1}{32} \sum_{i,j=0}^{3} \sum_{k=0}^{1} \chi(bba^{-i}c^{-j}d^{-k}a^i c^j d^k) \\
&= \frac{1}{8} \sum_{j=0}^{3} \sum_{k=1}^{2} \chi(c^{-j}d^{-k}c^j d^k) = \frac{1}{8} \left( \sum_{j=0}^{3} \chi_V(c^{-j}c^j) + \sum_{j=0}^{3} \chi(c^{-j}c^{-j}) \right) \\
&= \frac{1}{8}\chi_V(4 + 2 + 2c^2) = \frac{1}{4}\chi(3 + c^2)
\end{aligned}$$

Now let $V$ be the two-dimensional representation $\rho$ of $C_a$ given by
$$\rho(a) = \begin{pmatrix} 1 & 0 \\ 0 & 1 \end{pmatrix} \quad \rho(c) = \begin{pmatrix} i & 0 \\ 0 & -i \end{pmatrix} \quad \rho(d) = \begin{pmatrix} 0 & 1 \\ -1 & 0 \end{pmatrix}.$$
Then $\nu(\chi) = -1$ but $\nu_H(\hat{\chi}) = 1$.

## 5. The local indicator

We will see in this section that in order to compute the indicator in $H$, we only have to induce the $L_g$-module $V$ up to $\tilde{L}_g$, and not all the way up to $L$. We fix the notation $\tilde{V} := k\tilde{L}_g \otimes_{kL_g} V = \text{Ind}_{L_g}^{\tilde{L}_g} V$, and let $\hat{V} = \text{Ind}_{L_g}^{L} V$ as in Corollary 3.5. Also let $\Lambda_g = \frac{1}{|L_g|} \sum_{z \in L_g} z$ be the integral in $kL_g$ and let $\tilde{\Lambda}_g = \frac{1}{|\tilde{L}_g|} \sum_{z \in \tilde{L}_g} z$ be the integral in $k\tilde{L}_g$. For any $h \in H$, we write $h^{[2]} = \sum h_1 h_2$.

**Theorem 5.1.** Let $H$ be cocentral abelian as in (2.1), and assume that $\sigma$ is trivial. Choose $g \in G$, let $V$ be a simple $kL_g$-module and let $\tilde{V}, \hat{V}$ be the induced modules as above, with characters $\chi, \tilde{\chi}$, and $\hat{\chi}$ respectively.

1. If $g^2 = 1$ then
$$\nu_H(\hat{\chi}) = (\tau_{g,g^{-1}} * \chi^{[2]})(\Lambda_g).$$

2. If $o(g) \geq 3$ then
$$\nu_H(\hat{\chi}) = (\tau_{g,g^{-1}} * \tilde{\chi}^{[2]})(\tilde{\Lambda}_g) - (\tau_{g,g^{-1}} * \chi^{[2]})(\Lambda_g).$$

*Proof.* Let $\varphi = \tau_{g,g^{-1}}$; note $\varphi$ is a multiplicative character on $k\tilde{L}_g$ by Lemma 4.6.

1. If $g^2 = 1$ then $L_g = \tilde{L}_g = L_{g,g^{-1}}$ and by Proposition 4.12
$$\begin{aligned}\nu_H(\hat{\chi}) &= \frac{1}{|L_g|} \sum_{z \in L_g} \tau_{g,g^{-1}}(z)\chi(z^2) = \frac{1}{|L_g|} \sum_{z \in L_g} \varphi(z)(\chi^{[2]}(z)) \\ &= (\varphi * \chi^{[2]})(\frac{1}{|L_g|} \sum_{z \in L_g} z) = (\varphi * \chi^{[2]})(\Lambda_g).\end{aligned}$$

2. If $o(g) \geq 3$, there are two cases. First, assume that $g$ is $L$-non-real; that is, $L_{g,g^{-1}} = \emptyset$. Thus
$$\nu_H(\hat{\chi}) = 0 = (\varphi * \tilde{\chi}^{[2]})(\tilde{\Lambda}_g) - (\varphi * \chi^{[2]})(\Lambda_g).$$

Now assume $g$ is $L$-real; thus $L_g \neq \tilde{L}_g$. Then by Corollary 4.12
$$\begin{aligned}\nu_H(\hat{\chi}) &= \frac{1}{|L_g|} \sum_{z \in L_{g,g^{-1}}} \tau_{g,g^{-1}}(z)\chi(z^2) \\ &= \frac{1}{|L_g|} \sum_{z \in \tilde{L}_g} \tau_{g,g^{-1}}(z)\chi(z^2) - \frac{1}{|L_g|} \sum_{z \in L_g} \tau_{g,g^{-1}}(z)\chi(z^2).\end{aligned}$$

As in the case when $g^2 = 1$,
$$\frac{1}{|L_g|} \sum_{z \in L_g} \tau_{g,g^{-1}}(z)\chi(z^2) = (\varphi * \chi^{[2]})(\Lambda_g)$$

and
$$\frac{1}{|\tilde{L}_g|} \sum_{z \in \tilde{L}_g} \tau_{g,g^{-1}}(z)\tilde{\chi}(z^2) = (\varphi * \tilde{\chi}^{[2]})(\tilde{\Lambda}_g).$$

Therefore to complete the proof we need to show
$$\frac{1}{|\tilde{L}_g|} \sum_{z \in \tilde{L}_g} \varphi(z)\tilde{\chi}(z^2) = \frac{1}{|L_g|} \sum_{z \in \tilde{L}_g} \varphi(z)\chi(z^2).$$

Assume first that $\tilde{\chi}$ is a simple character of $\tilde{L}_g$. In this case the restriction of $\tilde{\chi}$ to $L_g$ is the sum of two characters $\chi$ and $\chi'$, where $\chi'$ is a conjugate of $\chi$. More precisely, we have $\chi' = \chi^d$ for (any) $d \in \tilde{L}_g \setminus L_g$, and for $y \in L_g$ we have $\chi'(y) = \chi^d(y) = \chi(dyd^{-1})$. In particular, if $z \in \tilde{L}_g$ then
$$\chi'(z^2) = \chi((dzd^{-1})^2).$$

Consequently

$$\begin{aligned}
\frac{1}{|\tilde{L}_g|} \sum_{z \in \tilde{L}_g} \varphi(z)\tilde{\chi}(z^2) &= \frac{1}{2|L_g|} \sum_{z \in \tilde{L}_g} \varphi(z)(\chi(z^2) + \chi((dzd^{-1})^2)) \\
&= \frac{1}{2|L_g|} \sum_{z \in \tilde{L}_g} (\varphi(z)\chi(z^2) + \varphi(dzd^{-1})\chi((dzd^{-1})^2)) \\
&= \frac{1}{|L_g|} \sum_{z \in \tilde{L}_g} \varphi(z)\chi(z^2).
\end{aligned}$$

If $\tilde{\chi}$ is not simple, then it is the sum of two characters $\xi$ and $\xi'$ of $\tilde{L}_g$, each of which is equal to $\chi$ when restricted to $L_g$. In particular $\xi(z^2) = \xi'(z^2) = \chi(z^2)$ for $z \in \tilde{L}_g$, and we find that

$$\frac{1}{|\tilde{L}_g|} \sum_{z \in \tilde{L}_g} \varphi(z)\tilde{\chi}(z^2) = \frac{1}{2|L_g|} \sum_{z \in \tilde{L}_g} \varphi(z)(\xi(z^2) + \xi'(z^2)) = \frac{1}{|L_g|} \sum_{z \in \tilde{L}_g} \varphi(z)\chi(z^2)$$

as desired. $\square$

In order to restate our theorem in a way analogous to the statement of the Frobenius-Schur theorem in [Se, Ch 13], we recall that a bilinear form may be defined on $H^*$ as follows: choose $\Lambda \in \int_H$ with $\varepsilon(\Lambda) = 1$. Then for $\phi, \psi \in H^*$, let

$$(\phi|\psi)_H := (S\phi * \psi)(\Lambda) = \sum S\phi(\Lambda_1)\psi(\Lambda_2).$$

The corollary is an immediate consequence of this definition and our main theorem.

**Corollary 5.2.** *Let $H$ be as in Theorem 5.1.*
1. *If $g^2 = 1$ then*
$$\nu_H(\hat{\chi}) = (\tau_{g,g^{-1}}^{-1}|\chi^{[2]})_{L_g}.$$
2. *If $o(g) \geq 3$ then*
$$\nu_H(\hat{\chi}) = (\tau_{g,g^{-1}}^{-1}|\tilde{\chi}^{[2]})_{\tilde{L}_g} - (\tau_{g,g^{-1}}^{-1}|\chi^{[2]})_{L_g}.$$

**Remark 5.3.** The theorem takes a simpler form when we make assumptions on $\tau$ as well. For example, assume that $\tau_{g,g^{-1}} = \psi^2$, for some $\psi \in \mathcal{G}((k\tilde{L}_g)^*)$ (in particular, this will be true if $\tau$ is trivial, or if $|L|$ is odd, as in Theorem 4.13).

Then since $(kL_g)^*$ and $(k\tilde{L}_g)^*$ are commutative,

$$(\tau_{g,g^{-1}} * \chi^{[2]})(\Lambda) = (\psi^{[2]} * \chi^{[2]})(\Lambda) = (\psi * \chi)^{[2]}(\Lambda) = (\psi * \chi)(\Lambda^{[2]}) = \nu(\psi * \chi).$$

Similarly $(\tau_{g,g^{-1}}^{-1}|\tilde{\chi}^{[2]})_{\tilde{L}_g} = \tilde{\nu}(\psi * \tilde{\chi})$. Thus we obtain a difference of ordinary indicators.

We now define the "local indicator" of an $L_g$-module $V$.

**Definition 5.4.** Fix $g \in G$ and let $V$ be a $L_g$-module. As before let $\tilde{V} := Ind_{L_g}^{\tilde{L}_g}(V)$, let $\chi$ and $\tilde{\chi}$ be the characters of $V$ and $\tilde{V}$, and let $\tilde{\nu}$ and $\nu$ be the classical Frobenius-Schur indicators of $\chi$ and $\tilde{\chi}$ for the group algebras $kL_g$ and $k\tilde{L}_g$. We define the *local indicator* $\nu_g(\chi)$ as follows:

(i) If $g^2 = 1$, we define $\nu_g(\chi) := \nu(\chi)$.
(ii) If $o(g) \geq 3$, we define $\nu_g(\chi) := \tilde{\nu}(\tilde{\chi}) - \nu(\chi)$.

Note that $\nu_g$ is identically zero if $g$ is $L$-non–real.

We can now restate our main theorem in terms of the local indicator, using Remark 5.3.

**Theorem 5.5.** *Let $H$ be a cocentral abelian extension with trivial cocycles. Let $V$ be a simple $L_g$-module, for some $g \in G$, and let $\hat{V}$ be the module obtained by inducing $V$ to $L$, so that $\hat{V}$ is an $H$-module. As before let $\chi$ and $\hat{\chi}$ be the characters of $V$ and $\hat{V}$. Then*
$$\nu_H(\hat{\chi}) = \nu_g(\chi).$$
*That is, the $H$-indicator of $\hat{\chi}$ is just the local indicator of $\chi$.*

The theorem shows that, using Corollary 3.4, the local indicator fits into the following commutative triangle of additive abelian groups:

(5.6)
$$\begin{array}{ccc} K_0(kL_g) & \xrightarrow{\phi} & K_0(H(\mathcal{O})) \\ {\scriptstyle \nu_g} \searrow & & \swarrow {\scriptstyle \nu} \\ & \mathbb{Z} & \end{array}$$

where $\mathcal{O}$ is the $L$-orbit of $g$, $K_0(H(\mathcal{O}))$ is the subgroup of $K_0(H)$ generated by the simple $H(\mathcal{O})$-modules, and $\phi$ is as in 3.4.

## 6. A criterion for positivity

This section is motivated by the well-known fact that for $G = S_n$, the symmetric group on $n$ letters, every simple character is real. We wish to extend this to the Drinfeld double, as well as consider more generally when our abelian extensions have this property.

**Theorem 6.1.** *Let $H = (kG)^* \# kL$ be a cocentral abelian extension as above, with trivial cocycle and dual cocycle. Assume that the following three conditions hold:*
 (i) $\tilde{\nu}(V) = 1$ *for all simple modules $V$ of $\tilde{L}_g$, for all $g \in G$, and*
 (ii) *If $d \in \tilde{L}_g \backslash L_g$, then the (conjugation) action of $d$ on $L_g$-modules sends $V$ to $V^*$.*
 (iii) *If $o(g) \geq 3$, then $\tilde{L}_g \neq L_g$.*
*Then every simple $H$-module has indicator equal to 1.*

*Proof.* Take a simple $L_g$-module $V$ with character $\zeta$. According to Theorem 5.5, we must show that $\hat{\nu}(V) = 1$ where $\hat{\nu}(V)$ is as in Definition 5.4. Set $\omega = \tilde{\zeta} = Ind_{L_g}^{\tilde{L}_g}(\zeta)$.

Suppose first that $g^2 = 1$. Then $L_g = \tilde{L}_g$, whence $\nu(V) = 1$ by (i). As $\hat{\nu}(V) = \nu(V)$ in this case, we are done by Theorem 5.5.

Now assume that $g$ has order at least 3. By assumption ($iii$), we know that $\tilde{L}_g \neq L_g$. The first case is when $V^* \neq V$. Then $\nu(\zeta) = 0$ and $\zeta^u$ is distinct from $\zeta$ by assumption ($ii$). The latter fact implies that $\omega$ is simple by elementary character theory, so that $\tilde{\nu}(\omega) = 1$ by assumption ($i$). Then again by Theorem 5.5, $\hat{\nu}(V) = \tilde{\nu}(\omega) - \nu(\zeta) = 1 - 0 = 1$ as required.

The last case is that $V^* = V$. This time elementary character theory and (ii) tell us that $\omega$ is the sum of a pair of simple characters $\tau_1$ and $\tau_2$, say, and we get $\tilde{\nu}(\tau_1) = \tilde{\nu}(\tau_2) = 1$ by (i) once more. Furthermore $\zeta = \operatorname{Res}_{L_g}^{\tilde{L}_g}(\tau_i)$ and so $\nu(\zeta) = \tilde{\nu}(\tau_i)$. Thus, $\hat{\nu}(V) = \tilde{\nu}(\omega) - \nu(\zeta) = 2 - 1 = 1$. □

When $H = D(G)$, we do not need condition (iii).

**Corollary 6.2.** *Let $H = D(G)$ and assume that the following two conditions hold:*
*(i) $\tilde{\nu}(V) = 1$ for all simple modules $V$ of $\tilde{C}_g$, for all $g \in G$, and*
*(ii) If $d \in \tilde{C}_g \backslash C_g$, then the (conjugation) action of $d$ on $C_g$-modules sends $V$ to $V^*$.*
*Then every simple $H$-module has indicator equal to 1.*

*Proof.* We only need to show that (iii) is automatically satisfied when $G$ acts on itself by conjugation. However this follows by a classical result of Burnside: for any $g \in G$, $C_{g,g^{-1}} \neq \emptyset$ if and only if all simple characters of $G$ are real valued. But by (i) when $g = 1$, we know that all simple characters of $G$ are real-valued. Thus $C_{g,g^{-1}} \neq \emptyset$ and so (iii) holds. □

Using the criteria in Corollary 6.2, we prove the result about $S_n$ mentioned in the introduction.

**Theorem 6.3.** *Every simple $D(S_n)$-module has indicator equal to 1.*

Every element $x \in S_n$ is the product of disjoint cycles which we indicate symbolically by

$$(6.4) \qquad x = \prod_{1 \leq i \leq n} t_i^{e_i}$$

which means that $x$ is the product of $e_i$ cycles of length $t_i$, $1 \leq i \leq n$. Of course, each $e_i$ is a non-negative integer. It is well-known that we then have

$$(6.5) \qquad C_x \cong \prod_{1 \leq i \leq n} \mathbb{Z}_{t_i} \wr S_{e_i}.$$

Thus $C_x$ is isomorphic to the direct product of (regular) wreathed products $\mathbb{Z}_{t_i} \wr S_{e_i}$, the latter being a semi-direct product of a symmetric group $S_{e_i}$ regularly permuting a basis of elements of a homogeneous abelian group isomorphic to $\mathbb{Z}_{t_i}^{e_i}$. We may thus represent $C_x$ as a semidirect product $C_x = A \rtimes P$ where $A$ is an abelian normal subgroup isomorphic to the direct product of the $\mathbb{Z}_{t_i}^{e_i}$ and $P$ is the direct product of the corresponding symmetric groups $S_{e_i}$. In this notation, if $x$ has order greater than two then the extended centralizer is a semidirect product

$$(6.6) \qquad \tilde{C}_x \cong A \rtimes (\langle u \rangle \times P)$$

where $u$ is an involution that inverts the elements of $A$ and commutes with $P$.

We can now show property (i).

**Lemma 6.7.** *Each simple $\tilde{C}_x$-module $V$ satisfies $\nu(V) = 1$.*

*Proof.* If $x = 1$ then $\tilde{C}_x = S_n$ and the result is well-known. So we may assume that $x$ has order at least 2. We make use of the standard description of simple characters of groups which are split extensions of a subgroup by a normal abelian subgroup (see [CR, Proposition 11.8]). Namely, pick a simple character $\xi$ of the normal abelian subgroup $A$, let $X$ be a complement to $A$ in $\tilde{C}_x$, and let $T = X_\xi$ be the stabilizer of $\xi$ in $X$. We may take $X = \langle u \rangle \times P$ or $P$ according as $x$ has order greater than 2 or not. Extend $\xi$ to a linear character of $AT$, denoted by $\xi_1$, by setting $\xi(t) = 1$ for $t \in T$. Pick a simple character $\psi$ of $T$, and inflate it to a simple character of $AT$, also denoted by $\psi$. Then the induced character $\zeta = Ind_{AT}^{\tilde{C}_x}(\xi_1 \psi)$ is a simple character of $\tilde{C}_x$, and every simple character of $\tilde{C}_x$ arises in this manner. An element $g \in \tilde{C}_x$ is uniquely expressible in the form $g = as$ with $a \in A$ and $s \in X$. Then $g^2 = (asas^{-1})s^2$ and we have

$$(6.8) \qquad \tilde{\nu}(\zeta) = |AT|^{-1} \sum_{s \in X} \psi(s^2) \sum_{a \in A} \xi \cdot \xi^s(a)$$

This holds since

$$\begin{aligned}
\tilde{\nu}(\zeta) &= |\tilde{C}_x|^{-1} \sum_{g \in \tilde{C}_x} \zeta(g^2) \\
&= |\tilde{C}_x|^{-1} |AT|^{-1} \sum_{g \in \tilde{C}_x, h \in \tilde{C}_x} \xi_1 \psi(hg^2 h^{-1}) \\
&= |AT|^{-1} \sum_{g \in \tilde{C}_x} \xi_1 \psi(g^2) \\
&= |AT|^{-1} \sum_{a \in A, s \in X} \xi(asas^{-1}) \psi(s^2) \\
&= |AT|^{-1} \sum_{s \in X} \psi(s^2) \sum_{a \in A} \xi \cdot \xi^s(a).
\end{aligned}$$

where $\xi^s$ is the $s$-conjugate of $\xi$, that is, $\xi^s(a) = \xi(sas^{-1})$ for $a \in A$. By orthogonality of characters, the inner sum in (6.8) is 0 unless $\xi^s = \xi^*$, the dual character, in which case it is $|A|$.

Suppose first that $\xi$ is real. Then $\xi^s = \xi^*$ if, and only if, $s$ lies in $T$. Equation (6.8) then implies that

$$(6.9) \qquad \tilde{\nu}(\zeta) = \nu(\psi).$$

On the other hand, if $\xi$ is not real then $\xi^u = \xi^*$ is distinct from $\xi$, so $u$ does not lie in $T$ and the set of elements $s \in S$ satisfying $\xi^s = \xi^*$ is precisely the coset $uT$. Noting that $(ut)^2 = t^2$ for $t \in T$, we see from (6.8) that (6.9) continues to hold. So in fact (6.9) holds in all cases. But it is easy to see that $T$ is isomorphic to either a direct product of symmetric groups (if $\xi$ is not real), or a direct product of $\langle u \rangle$ and some symmetric groups. In either case we get $\nu(\psi) = 1$, and the lemma now follows from (6.9). □

We now show that (ii) holds for $S_n$, finishing the proof of Theorem 6.3.

**Lemma 6.10.** *Assume that $o(g) \geq 3$. Then the (conjugation) action of $u$ on $C_g$-modules $V$ coincides with the map sending $V$ to its contragredient $V^*$.*

*Proof.* Set $C = C_g$. First some reductions: we may assume that $V$ is simple with character $\zeta$, and try to prove that

(6.11) $$\zeta^u = \zeta^*.$$

As $u$ normalizes each of the direct factors of $C$ (cf (6.5)) we may also assume without loss that $C$ coincides with one of these factors. Thus, to simplify notation, we may take $A = (\mathbb{Z}_k)^n$ and $P = S_n$.

Once again the simple characters of $C$ are constructed in the same manner as described in the proof of Lemma 6.10. Thus we have

(6.12) $$\zeta = Ind_{AT}^C(\xi_1 \cdot \psi)$$

where $\xi$ is a simple character of $A$ with stabilizer $T$ in $P$ and extended to a character of $AT$ by assigning the value 1 at elements of $T$, and $\psi$ is a simple character of $T$. We have $T \cong S_m$ for some $m$. Now because $u$ inverts $A$ then $\xi^u = \xi^*$; since $u$ commutes with $T$ then $\xi_1^u = \xi_1^*$ and $\psi^u = \psi = \psi^*$, the latter because all characters of symmetric groups are real-valued. Thus (6.12) yields

$$\zeta^u = Ind_{AT}^C(\xi_1^u \cdot \psi^u) = Ind_{AT}^C(\xi_1^* \cdot \psi^*) = Ind_{AT}^C(\xi_1 \cdot \psi)^* = \zeta^*$$

thereby establishing (6.11). □

We give another application of Theorem 6.1.

**Lemma 6.13.** *Let $G = D_{2n}$ be the dihedral group of order $2n$ and let $H = D(G)$. Then $\nu_H(\hat{\chi}) = 1$ for all simple characters $\hat{\chi}$ of $H$.*

*Proof.* Again it is known that $\nu(\chi) = 1$ for all irreducible characters of $D_{2n}$. We need to verify the two conditions of Corollary 6.2.

Let $H$ be the cyclic subgroup of $G$ of order $n$. If $h \in H$, then $h$ is inverted under conjugation, and $\tilde{C}_h = G$, so certainly (i) holds. If $g \in G \backslash H$, then $g$ has order 2, and so $\tilde{C}_g = C_g$, which is isomorphic to $\mathbb{Z}_2$ or $(\mathbb{Z}_2)^2$ depending on whether $n$ is odd or even. In either case again (i) is satisfied.

For (ii) we only need to consider elements $h$ in $H$. If $[\tilde{C}_h : C_h] = 2$ then $C_h = H$ and so if $g \in G \backslash H$, $g$ conjugates each conjugacy class of $C_h$ to its inverse class. But this is equivalent to (ii). □

We remark that we could also have proved Lemma 6.13 directly from Theorem 4.4, since the representations of $D_{2n}$ are well-known.

In fact the proof of Lemma 6.13 applies to a generalized dihedral group, that is, a group of the form $G = \langle H, t \rangle$, where $H$ has index 2 in $G$ and $tht = h^{-1}$ for all $h \in H$. Combining this with Theorem 6.3, we have shown:

**Corollary 6.14.** *Let $G$ be any direct product of groups of the form generalized dihedral, elementary 2-group, or symmetric group. Then $\nu(\chi) = 1$ for all simple characters of $D(G)$.*

**Remark 6.15.** It is likely that Theorem 6.3 remains true for other classes of groups, in particular for the Weyl groups associated to semisimple root systems. Theorem 6.3 and Corollary 6.14 provide the proof if the simple components are of type $A$. It should not be difficult to extend our proof to the other types, in particular for types $B$ and $D$ where the Weyl group is just an extension of a symmetric group by an elementary abelian 2-group. It is well-known, for example, that one ingredient of our proof - that the rationals are a splitting field for $S_n$ - holds for the other Weyl groups (see [GP]).

## 7. Hopf algebras of dimension 16 revisited

The semisimple Hopf algebras of dim 16 have recently been classified, in [K]: there are exactly 16 non-trivial such Hopf algebras. All of them can be described as cocentral abelian extensions with $L = \mathbb{Z}_2 = \langle z \rangle$, and thus the results in this paper apply. A major tool in [K] was the computation of the rings $K_0(H)$. We show here that using the Schur indicator provides a somewhat simpler invariant for the classification.

We first need a few consequences of Theorem 5.5. As before, $V$ denotes a simple $L_g$-module with character $\chi$.

**Lemma 7.1.** *Assume that $\sigma$ is trivial, $g^2 = 1$, and $L_g$ is abelian. Then $\nu_H(\hat{\chi})$ is nonnegative.*

*Consequently if $G = k(\mathbb{Z}_2)^n$ and $L$ is a cyclic group, then for any $\hat{\chi} \in Irr(H)$, $\nu_H(\hat{\chi})$ is nonnegative.*

*Proof.* Since $\tilde{L}_g = L_g$ is abelian, $V$ is one-dimensional and so $\chi$ is a group-like element. Thus $\chi^{[2]} = \chi^2$ and therefore by Theorem 5.1

$$\nu_H(\hat{\chi}) = (\tau_{g,g^{-1}} * \chi^2)(\Lambda_g) = \begin{cases} 1 & \text{if } \chi^2 = \tau_{g,g^{-1}}^{-1} \\ 0 & \text{otherwise} \end{cases}$$

When $L$ is cyclic, it is known that $H^2(L, (k^G)^\times) = 0$. Thus as in [Na, 1.2.6] we may assume that $\sigma$ is trivial, and so $\nu_H$ is non-negative by the first argument. □

**Lemma 7.2.** *Assume that $L$ is an abelian group of exponent $2n$, where $n$ is odd, and that $\sigma$ and $\tau$ are trivial. Then for any $\hat{\chi} \in Irr(H)$, $\nu_H(\hat{\chi})$ is nonnegative.*

*Proof.* Assume as before that $\hat{\chi}$ is the character of $\hat{V}$, for $V$ an irreducible $L_g$-module, for some $g \in G$. Since $L$ is abelian, $V$ is one-dimensional and therefore $\nu(\chi) = 0$ or 1.

If $g^2 = 1$, then $\nu_H(\hat{\chi}) = \nu(\chi)$ is nonnegative by Theorem 5.5.

If $o(g) \geq 3$ and $L_g = \tilde{L}_g$ then $\nu_H(\hat{\chi}) = 0$.

If $o(g) \geq 3$ and $L_g \neq \tilde{L}_g$ then $\tilde{\chi} = \xi + \xi'$ where $\xi$ and $\xi'$ are one-dimensional characters of $\tilde{L}_g$ since $\tilde{L}_g$ is abelian. Therefore $\tilde{\nu}(\tilde{\chi}) = \nu(\xi) + \nu(\xi') \geq 0$.

If $\nu(\chi) = 0$ then $\nu_H(\hat{\chi}) = \tilde{\nu}(\tilde{\chi}) - \nu(\chi) \geq 0$. If $\nu(\chi) = 1$ then $\chi^2 = \varepsilon_{kL_g}$. Since $\exp(L) = 2n$, $n$ odd, this implies that $\xi^2 = (\xi')^2 = \varepsilon_{k\tilde{L}_g}$. Therefore $\nu_H(\hat{\chi}) = \tilde{\nu}(\tilde{\chi}) - \nu(\chi) = 2 - 1 = 1$. □

**Corollary 7.3.** *Let $H$ be a semisimple Hopf algebra of dimension $16$. Then $\nu(\chi)$ is known for all simple $H$-modules. Moreover, all such $H$ with $\mathcal{G}(H) \not\cong \mathbb{Z}_2 \times \mathbb{Z}_2 \times \mathbb{Z}_2$ can be distinguished by $\nu_H$, $\mathcal{G}(H)$, $\mathcal{G}(H^*)$ and, in the case of $\mathcal{G}(H) \cong \mathbb{Z}_4 \times \mathbb{Z}_2$, the action $\rightharpoonup$ of $L$ on $G$.*

*Proof.* We consider the possible groups of group-like elements. Following Table 1 of [K], there are four Hopf algebras with $\mathcal{G}(H) \cong \mathbb{Z}_2 \times \mathbb{Z}_2 \times \mathbb{Z}_2$, seven with $\mathcal{G}(H) \cong \mathbb{Z}_4 \times \mathbb{Z}_2$, two with $\mathcal{G}(H) \cong D_8$, and three with $\mathcal{G}(H) \cong \mathbb{Z}_2 \times \mathbb{Z}_2$. If $\mathcal{G}(H)$ is abelian of order 8, we may use $K = (k\mathcal{G}(H))^*$, and if $\mathcal{G}(H)$ is abelian of order 4 we use $K = (kD_8)^*$. If $\mathcal{G}(H) \cong D_8$, then in the two possible examples, one has $K = (kD_8)^*$ and the other has $K = (kQ_8)^*$.

If $\mathcal{G}(H) \cong \mathbb{Z}_2 \times \mathbb{Z}_2 \times \mathbb{Z}_2$ then $H$ has two simple degree 2 characters, $\chi_1$ and $\chi_2$, and they are self-dual (see [K, Section 3.2]). Therefore by Lemma 7.1, $\nu(\chi_1) = \nu(\chi_2) = 1$.

If $\mathcal{G}(H) \cong \mathbb{Z}_4 \times \mathbb{Z}_2 = \langle g_1 \rangle \times \langle g_2 \rangle$ then there are three different types of actions: $z \rightharpoonup_a e_{i,j} = e_{i+2j,j}$, $z \rightharpoonup_b e_{i,j} = e_{-i,j}$ and $z \rightharpoonup_c e_{i,j} = e_{i,j+i}$, where $\{e_{i,j}\}_{i=0,1,2,3;j=0,1}$ is the basis dual to $\{g_1^i g_2^j\}$. $H$ has two simple degree 2 characters $\chi_1$ and $\chi_2$ (see [K, Section 3.2]). Using Lemma 4.1 we can compute that

$$\Lambda_a^{[2]} = \begin{cases} e_{0,0} + e_{2,0} + \frac{1}{2}(e_{0,1} + e_{2,1}) - \frac{1}{2}(e_{1,1} + e_{3,1}) & \text{if } H = H_{a:1} \\ e_{0,0} + e_{2,0} + \frac{1}{2}(e_{0,1} + e_{2,1}) + \frac{1}{2}(e_{1,1} + e_{3,1}) & \text{if } H = H_{a:y} \\ e_{0,0} + e_{2,0} + e_{0,1} + e_{2,1} + \frac{1}{2}(e_{1,0} + e_{3,0}) - \frac{1}{2}(e_{1,1} + e_{3,1}) & \text{if } H = H_{b:1} \\ e_{0,0} + e_{2,0} + \frac{1}{2}(e_{1,0} + e_{3,0}) + \frac{1}{2}(e_{1,1} + e_{3,1}) & \text{if } H = H_{b:y} \\ e_{0,0} + e_{2,0} - \frac{1}{2}(e_{1,0} + e_{3,0}) - \frac{1}{2}(e_{1,1} + e_{3,1}) & \text{if } H = H_{b:x^2 y} \end{cases}$$

Then $\nu_{H_{a:1}}(\chi_1) = -1$ and $\nu_{H_{a:1}}(\chi_2) = 1$; $\nu_{H_{a:y}}(\chi_1) = \nu_{H_{a:y}}(\chi_2) = 1$; $\nu_{H_{b:1}}(\chi_1) = 1$ and $\nu_{H_{b:1}}(\chi_2) = -1$; $\nu_{H_{b:y}}(\chi_1) = \nu_{H_{b:y}}(\chi_2) = 1$; $\nu_{H_{b:x^2 y}}(\chi_1) = \nu_{H_{b:x^2 y}}(\chi_2) = -1$. If $H = H_{c:\sigma_0}$ or $H_{c:\sigma_1}$, the degree 2 simple characters are not self-dual and therefore $\nu_H(\chi_1) = \nu_H(\chi_2) = 0$. Here $\mathcal{G}(H_{b:1}^*) = \mathcal{G}(H_{c:\sigma_1}^*) = \mathbb{Z}_2 \times \mathbb{Z}_2 \times \mathbb{Z}_2$; in the other five cases $\mathcal{G}(H^*) = \mathbb{Z}_4 \times \mathbb{Z}_2$.

If $\mathcal{G}(H) \cong D_8$ then $H = H_{C:1}$ or $H_E$ and $H$ can be written as an abelian extension as in Lemma 7.2 (see [K, Sections 3.3 and 3.4]). $H$ has three simple degree 2 characters, $\chi_1$, $\chi_2$ and $\chi_3$. If $H = H_{C:1}$ they are self-dual. Thus $\nu(\chi_1) = \nu(\chi_2) = \nu(\chi_3) = 1$. If $H = H_E$ then $\chi_2$ is selfdual and $\chi_1$ and $\chi_3$ are not. Thus $\nu(\chi_1) = \nu(\chi_3) = 0$ and $\nu(\chi_2) = 1$. In both of these examples, $\mathcal{G}(H^*) = \mathbb{Z}_2 \times \mathbb{Z}_2$.

If $\mathcal{G}(H) \cong \mathbb{Z}_2 \times \mathbb{Z}_2$ then $H = H_{B:1}, H_{B:X}$ or $H_{C:\sigma_1}$. In the case of $H = H_B$, $H$ has two simple degree 2 characters $\chi_1$ and $\chi_2$. By Lemma 4.1 $\nu_{H_{B:1}}(\chi_1) = \nu_{H_{B:1}}(\chi_2) = 1$; also $\nu_{H_{B:X}}(\chi_1) = -1$ and $\nu_{H_{B:X}}(\chi_2) = 1$. In both of these examples, $\mathcal{G}(H^*) \cong D_8$.

In the case of $H = H_C$, $H$ has three simple degree 2 characters $\chi_1$, $\chi_2$, and $\chi_3$. By Lemma 4.1, $\nu_{H_{C:\sigma_1}}(\chi_1) = \nu_{H_{C:\sigma_1}}(\chi_2) = \nu_{H_{C:\sigma_1}}(\chi_3) = 1$. In this last case, $\mathcal{G}(H^*) \cong \mathbb{Z}_2 \times \mathbb{Z}_2$. □

We give a different, more direct argument for $H_{B:X}$ in the next example.

**Example 7.4.** Let $H = H_{B:X} = kQ_8 \#^\alpha k\mathbb{Z}_2$, a smash coproduct of $kQ_8$ and $k\mathbb{Z}_2$ (see [Mo, 10.6.1] or [Ni, Proposition 3.8]), where $Q_8 = \langle a, b : a^4 = 1, b^2 = a^2, ba =$ $a^{-1}b\rangle$, $\mathbb{Z}_2 = \{1, g\}$. $H$ has the algebra structure of $kQ_8 \otimes k\mathbb{Z}_2$ and the following comultiplication, antipode and counit:

$$\Delta(x\#\delta_{g^k}) = \sum_{r+t=k} (x_1\#\delta_{g^r}) \otimes (\alpha_{g^r}(x_2)\#\delta_{g^t})$$
$$S(x\#\delta_{g^k}) = \alpha_{g^k}(S(x))\#\delta_{g^{-k}}$$
$$\varepsilon(x\#\delta_{g^k}) = \varepsilon(x)\delta_{g^k,1}$$

where $\delta_1 = \frac{1}{2}(1+g)$, $\delta_g = \frac{1}{2}(1-g)$, $x \in kQ_8$ and $\alpha : \mathbb{Z}_2 \to Aut(kQ_8)$ is given by $\alpha_1 = id$ and $\alpha_g$ is the automorphism interchanging $a$ and $b$; see [Ni, Erratum] and [K, Section 8].

The basis of $H$ is $\{a^{2l}a^ib^jg^k | i, j, k, l = 0, 1\}$. Then

$$\Lambda = \frac{1}{16} \sum_{i,j,k,l=0}^{1} a^{2l}a^ib^jg^k = \frac{1}{8} \sum_{i,j,l=0}^{1} a^{2l}a^ib^j\delta_1$$

$$\Delta\Lambda = \frac{1}{8} \sum_{i,j,l=0}^{1} \left(a^{2l}a^ib^j\delta_1 \otimes a^{2l}a^ib^j\delta_1 + a^{2l}a^ib^j\delta_g \otimes a^{2l}b^ia^j\delta_g\right)$$

$$\Lambda^{[2]} = \frac{1}{8} \sum_{i,j,l=0}^{1} \left(a^{4l}a^ib^ja^ib^j\delta_1 + a^{4l}a^ib^{i+j}a^j\delta_g\right)$$

$$= \frac{1}{4} \sum_{i,j=0}^{1} \left(a^ia^ia^{2ij}b^{2j}\delta_1 + a^ia^ja^{2(i+j)j}b^{i+j}\delta_g\right)$$

$$= \frac{1}{4} \sum_{i,j=0}^{1} \left(a^{2i+2j+2ij} + a^{(i+j)(2j+1)}b^{i+j}\right) = \frac{1}{4}\left(1 + 3a^2\right)\delta_1 + \frac{1}{4}\left(2 + ab + a^3b\right)\delta_g$$

Let $\pi_i$, $i = 1, 2$ be the simple $H$-modules of degree 2 (see [K, Section 8]). Then $\pi_i(a^2) = -I$, $\pi_1(\delta_1) = I$, $\pi_1(\delta_g) = 0$, $\pi_2(\delta_1) = 0$, $\pi_2(\delta_g) = I$. Thus $\pi_1\left(\Lambda^{[2]}\right) = \frac{1}{4}\pi_1\left(1 + 3a^2\right) = -\frac{1}{2}I$ and $\pi_2\left(\Lambda^{[2]}\right) = \frac{1}{4}\pi_2\left(2 + ab + a^3b\right)$. Therefore

$$\nu(\chi_1) = \chi_1\left(\Lambda^{[2]}\right) = -1$$
$$\nu(\chi_2) = \chi_2\left(\Lambda^{[2]}\right) = \frac{1}{4}(4 + 0 + 0) = 1$$

We compare this with the fact that for the usual group algebra $H = k[Q_8 \times \mathbb{Z}_2]$, both simple modules of degree two have $\nu(\chi) = -1$.

SYRACUSE UNIVERSITY, SYRACUSE, NY 13244-1150
*E-mail address*: ykashina@syr.edu

UNIVERSITY OF CALIFORNIA SANTA CRUZ, SANTA CRUZ CA 95064
*E-mail address*: gem@cats.ucsc.edu

UNIVERSITY OF SOUTHERN CALIFORNIA, LOS ANGELES, CA 90089-1113
*E-mail address*: smontgom@math.usc.edu